\documentstyle{amsppt}
\magnification=1200
\hcorrection{0.5in}
\hsize =14truecm
\baselineskip =18truept
\vcorrection{0.15in}
\vsize =20.5truecm

\nologo

\NoRunningHeads

\NoBlackBoxes

\topmatter
\document

\title {\bf On Projective varieties with nef
anticanonical divisors}
\endtitle
\author  {Qi  Zhang}    
\endauthor
\address{Department of mathematics,  University of
Missouri, Columbia, 
MO
65211. U.S.A}  \endaddress
\email{qi\@math.missouri.edu}\endemail
\endtopmatter

\heading{\bf $\S$1. Introduction}\endheading

The aim of this note is to prove a structure
theorem for projective varieties
with nef anticanonical divisors (the Main Theorem). In [18], we showed
that if $X$ is smooth and $-K_X$ is nef, then the
Albanese map
$\text{Alb}_X: X\rightarrow \text{Alb}(X)$ is surjective and has connected
fibers (i.e., it is a fiberspace map). In this note we  apply the  techniques which have been developed in
[2],[14] and [19] to prove the following;

\proclaim{Main Theorem} Let $X$ be a projective
variety and $D$  an effective $Q$-divisor on $X$  such
that the pair $(X,D)$
is  log canonical [10]  and $-(K_X+D)$  is nef.
Let $f:X-->Y$ be a dominant
rational map, where $Y$ is a smooth variety. Then either 
\roster
\item  $Y$ is  uniruled; or 
\item The Kodaira
dimension $\kappa (Y)=0$ 
 Moreover in this case, $f$ is semistable in
codimension 1 (see Definition 1).
\endroster
\endproclaim

{\it Remark 1}: In fact, we only need to assume that $(X,D)$ is log 
canonical along the general fibers of $f$. Also, in case (2), we  
 actually have $\kappa_{\sigma}(Y)=0$ [14].
 Moreover, our proof also shows that
in this case, every component of $D$ dominates $Y$. This fact was proved
by McKernan-Prokhorov in [20, Lemma 10.2] under the assumption of  the Log Minimal Model 
Conjecture.

{\it Remark 2}:  When $X$ is smooth and $f$ is a
morphism, the Main Theorem was essentially  
proved
in [18] by using the  relative deformation theory and
the mod $p$ reduction methods.
However, the proof we present here is  new and is
completely different in nature.

Applying the Main Theorem to the maximal rationally
connected fibration
of $X$ [3],[11] and  using a result  of Graber-Harris-Starr
[7], we   obtain the following:

\proclaim{Corollary 1}
Let $X$ be a projective variety and $D$  an effective
 $Q$-divisor  such that the pair $(X,D)$
is  log canonical and $-(K_X+D)$  is nef.
Then there exists a dominant rational map $f:X-->Y$  such that
\roster
\item $f$ is a fiberspace map which is also semistable in
codimension one.
\item The general fibers of $f$ are rationally
connected.
\item $Y$ is smooth with  $\kappa (Y)=0$.
\endroster
\endproclaim

{\it Remark 3}:  In the statement of  Corollary 1, I
believe that 
$Y$ should be birational to a variety $W$ with log
canonical singularities and $K_W$ being numerically
trivial.

\proclaim{Corollary 2}
Let $X$ be a projective variety and $D$  an effective
 $Q$-divisor  such that the pair $(X,D)$
is  log canonical and $-(K_X+D)$  is nef. Then the
Albanese map (from any smooth model of $X$) $\text{Alb}_X: X-->\text{Alb}(X)$
is a  fiberspace map which is also semistable in
codimension one.
\endproclaim

{\it Remark 4}: A conjecture of
Demailly-Peternell-Schneider [6] predicts that when $X$
is smooth and $-K_X$ is nef,  the Albanese map should
be a smooth morphism. When $\text{dim}X=3$,
 this has been verified by Peternell-Serrano [15]. Also,
Campana proved that when $X$ is a special variety,  its Albanese map
is multiple fibers free in codimension one [4].

\proclaim{Corollary 3}
Let $X$ be a projective variety and $D$  an effective
 $Q$-divisor  such that the pair $(X,D)$
is  log canonical and $-(K_X+D)$  is nef. Let $f:
X-->Y$ be the nef reduction of $-(K_X+D)$ [5],[16].
Then $Y$ is uniruled.
\endproclaim

{\it Remark 5}: The notion of special varieties was introduced and
 studied  by F. Campana in [4]. He also conjectured that compact
K\"{a}hler manifolds with $-K_X$ nef are special. S. Lu [12] proved the
 conjecture
for projective varieties. In particular, if $X$ is a
projective variety with $-K_X$ nef and if there is a surjective map $X\to Y$. 
Then $\kappa (Y)\leq 0$. Our  focus  however, is on the uniruledness of $Y$.

{\it Remark 6}: When $\text{dim}X\leq 4$ and $-K_X$ is
nef, we can easily derive some (rough) classification results
for those varieties from  Corollaries 1-3. A complete classification
 of such threefolds was given in [21].

The  tools we shall use in the proof of the Main
Theorem are the theory of  weak (semi) positivity of the direct images of
 (log) relative dualizing sheaves $f_{*}(K_{X/Y}+\Delta)$ (which has
been developed by  Fujita, Kawamata, Koll\'{a}r, Viehweg and others) plus 
the notion of the divisorial Zariski decomposition [1],[14].
Also, some  results in [2] are crucial to the proof.

\heading{\bf $\S$2. Proofs}\endheading

Before we start to prove  the Main Theorem and Corollaries 1-3, let us
first give the related definition.

\proclaim {Definition 1} Let $f: X-->Y$ be a dominant rational
map, where $X$ and $Y$ are normal proper varieties. We say that $f$
is semistable in codimension one if it satisfies the following condition
: 
\medskip

 Let   $\pi: Z\rightarrow X$ be a birational morphism which is an  
 {\it eliminating the indeterminacy of $f$}.    Let $P$ be a prime divisor  on $Y$. If $D$ is   
a non-reduced 
component of $g^{*}(P)$, where  $g=f\circ \pi$.  Then either $\text{codim}(g(D))\geq 2$ or $D$
is  $\pi$-exceptional.
\endproclaim

 We work over the complex number field $\Bbb C$ in  this paper.

\demo{Proof of the Main Theorem} The proof  is quite
similar to the proof of Proposition 1 in [19]. 

  Eliminating the indeterminacy of $f$ and
taking a log resolution, we have a smooth projective
variety $Z$, and the surjective morphisms
$g:Z\rightarrow Y$ and $\pi:Z\rightarrow X$, 

$$
\CD
Z  @>\pi>> X \\
@VVgV         @. \\
Y   @.     
\endCD
$$ such that $K_{Z}={\pi}^{*}(K_{X}+D)+\sum e_{i}E_{i}$ with
$e_{i}\geq -1$, where $\sum E_{i}$ is a divisor with 
normal crossing.

 Let $\delta L$ be an ample $Q$-divisor (with small
$\delta>0$) on $Z$.
 Since $-(K_X+D)$ is nef, $-{\pi}^{*}(K_X+D)+\delta L$ is 
ample. Let $H$ be an  ample $Q$-divisor on
$Y$ such that $A=-{\pi}^{*}(K_X+D)+\delta L-H$ is again
an ample $Q$-effective divisor. We may also assume 
that $\text{Supp}(A+\sum E_i)$ is a 
divisor with simple normal crossing. Let
$\Delta = A+\sum {\epsilon}_{i}E_{i}$, where ${\epsilon}_{i}E_{i}=
-e_{i}E_{i}$ if $e_{i}<0$  and  ${\epsilon}_{i}E_{i}=\{-e_{i}\}E_{i}$ if $e_{i}\geq 0$.
 We may assume that  the pair
 $(X, \Delta)$ is log canonical.  We have
$$K_{Z/Y}+\Delta \sim_{Q} \sum 
m_{i}E_{i}-g^{*}(K_Y+H)+\delta L$$ where $m_i={\lceil
e_{i} \rceil}$ (for $e_i\geq 0$) are non-negative
integers (where $\{,\}$ is the fractional part and
$\lceil, \rceil$ is the round up).

By the same arguments as in [19] (using the semistable
reduction theorem, the covering trick
and the flattening of $g$), there exist a finite
morphism
$p:Y'\rightarrow Y$  and the induced morphism
$g':Z'\rightarrow Y'$ from a desingularization 
$Z'\rightarrow Z\times_{Y}Y'$ such that 
 it is semistable over $Y'\setminus B$ with $\text{codim}(B)\geq
2$. Let $Z'\rightarrow Z$ be the induced morphism.

$$
\CD
Z'    @>q>> Z  @>\pi>> X \\
@Vg'VV   @VgVV   @. \\
Y'    @>p>>    Y @.
\endCD
$$

By the ramification formula and the fact that the relative dualizing
sheaves are stable under the base change. By [19] we have

$$K_{Z'/Y'}+\Delta'\sim_{Q} \sum
n_{k'}E'_{k'}-{g'}^{*}p^{*}(K_Y+H)+q^{*}(\delta
L)-V+G$$ where
\roster
\item  $ \sum n_{k'}E'_{k'}$ is Cartier and
$\text{Supp}(\sum E'_{k'})$ is $q\circ
\pi$-exceptional.
\item $V$ is an effective Cartier divisor which is
$g'$-vertical ($V$ contains the divisors from the non-reduced components 
of the fibers of $g$).
\item $G$ is also Cartier and $q\circ
\pi$-exceptional ($\text{codim}g'(G)\geq 2$).
\item If $E'_{k'}$ is  $g'$-horizontal, then
$n_{k'}\geq 0$ ($g'$-horizontal divisors are not the ramification divisors
of $q$). 
\item  $(Z',\Delta')$ is log canonical.
\item $p^{*}H=2H'$ where $H'$ is an  ample Cartier
divisor.
\item If a divisor $D$ on $Z'$ with
$\text{codim}(g'(D))\geq 2$, then $D$ is
 $q\circ \pi$-exceptional.
\endroster 

The following weak positivity theorem is a
generalization of the results due to Kawamata  [8], [9]
and Viehweg [17], and it was proved first by F. Campana [4, Theorem 4.13]
 (see also [12]).

\proclaim{Lemma 1} Let $U$ be a smooth projective
variety and $M$ be an effective divisor  such that the
pair $(U,M)$
is  log canonical. Let $\phi: U\rightarrow W$ be a
surjective morphism, where $W$ is smooth.
 Then $\phi_{*}(m(K_{U/W}+D))$ is torsion free and
weakly positive for any $m\in Z_{>0}$ such that
$m(K_{U/W}+D)$ is Cartier.
\endproclaim

Let us choose a positive integer $m>0$ such that $m(\delta
L)$ is Cartier.

Since all the
$g'$-horizontal divisors $E'_{k'}$ and
${q}^{*}(m\delta L)$ have non negative
coefficients, $$g'_{*}(m(
\sum
n_{k'}E''_{k'}-{g'}^{*}p^{*}(K_Y)-{g'}^{*}(2H')+{q}^{*}(\delta
L)-V+G))$$
is a non-zero sheaf. Let $$\omega =m(\sum
n_{k'}E'_{k'}-{g'}^{*}p^{*}(K_Y)+{q}^{*}(\delta L)-V+G)$$

By Lemma 1, $g'_{*}(\omega)\otimes \Cal O_{Y''}(-2mH')$
is torsion free and weakly positive.

Thus there exists an integer $n\in Z_{>0}$ such that
$$\Cal
{\hat {\Cal S}}^{n}{g'}_{*}(\omega )\otimes \Cal
O_{Y'}(-2mnH')\otimes \Cal O_{Y'}(nmH')$$   is 
 generically generated by
global sections  for some $n>0$,  where ${\hat {\Cal
S}^{n}}$ denotes 
the reflexive hull of ${\Cal S}^{n}$. This implies that there is a
divisor $B'$ on $Z'$ with  $\text{codim}g'(B')\geq 2$
 and an inclusion $$\Cal O_{Z''}\rightarrow
{\omega}^{n}\otimes {g'}^{*}(-mnH')+B'$$

Therefore

$$M=mn(\sum
n_{k'}E'_{k'}-{g'}^{*}p^{*}(K_X)-{g'}^{*}(H')+{q}^{*}(\delta
L)-V+G)+B'$$ is effective (see [19]).

On the other hand, we can  choose a  family of general
complete intersection curves $C$ on $Z'$ such that
$C$ does not intersect
with the exceptional locus of $Z'\rightarrow X$ (such
as ${E'}_{k'}$,$G$ and $B'$), e.g.,pull-back the
general
 complete intersection curves from $X$.

 Thus we have 
$C\cdot {g'}^{*}p^{*}(K_Y)\leq C
\cdot ({q}^{*}(\delta L)-V)$. i.e., there
exists a covering family of curves $C'$ on $Z$
such that $C'\cdot g^{*}K_Y\leq C'\cdot (\delta
L-q_{*}V)$. 

Let $\delta \rightarrow 0$, 
we get $C'\cdot g^{*}K_Y\leq 0$. If $C'\cdot
g^{*}K_Y<0$, $Y$ is uniruled [13]. Let us assume
that $C'\cdot  g^{*}K_Y=0$ (in particular, $(\pi\circ q)_{*}V=0$ and
$f$ is semistable in codim 1), and $Y$ is not
uniruled.

By a result in [2, Corollary 0.3], $K_Y$ is pseudo-effective.
Let $K_Y\equiv P+N$ be the divisorial Zariski
decomposition of $K_Y$, where $P$ is the positive part
(nef in codimension 1) and $N$ is the negative part
(an effective $R$-divisor) [1],[14]. Since the family
of the curves $g(C')$ also forms a {\it connecting}  family
[2] and $g(C')\cdot K_Y=0$, we have $p\equiv 0$ [2, Theorem 9.8].
Thus $\kappa_{\sigma}(Y)=0$  and hence $\kappa(Y)=0$ by a 
result of N. Nakayama [14, Proposition 6.2.8]. q.e.d.
\enddemo

\demo{Proof of  Corollary 1} Let $C$ be a general
complete intersection curve on $X$.
If $C\cdot (K_X+D)<0$, then  $X$ is uniruled  and there
exists a
non-trivial maximal rationally connected fibration
$f:X-->Y$ [3], [11].
 By a result in [7],  $Y$ is not uniruled if $\text{dim}Y>0$. If $Y$
is a point, then $X$ is rationally connected.
 If $\text{dim}Y>0$, then $\kappa(Y)=0$ by the Main
Theorem. 

So let us assume that $C\cdot (K_X+D)=0$ and $X$ is not uniruled. Thus we have $D=0$ and 
$K_X$ is numerically trivial. Let $\pi: Z\rightarrow X$ be a log resolution. Then $K_Z+\Delta={\pi}^{*}K_X+E$.
where $\Delta$ and $E$ are effective and  $\pi$-exceptional. Let $C'$ be the family of the curves from the pull-back 
of $C$. Then  $C'\cdot (K_Z+\Delta)=0$. Since $K_Z$ is pseudo-effective, we have $K_Z=P+N$ (the divisorial Zariski
decomposition of $K_Z$). Since $\text{Supp}(\Delta)$ is $\pi$-exceptional, $P$ is also the positive part of the divisorial Zariski
decomposition of $K_Z+\Delta$. On the other hand, $C$ forms a connecting family. Therefore we have $P\equiv 0$ and hence  
$K_{\sigma}(Z)=\kappa(Z)=0$ [14, Proposition 6.2.8]. q.e.d.
\enddemo

\demo{Proof of  Corollary 2} Since for any subvariety
$W$ of an abelian variety,  
$\kappa(W)\geq 0$ and $\kappa (W)=0$ if and only if $W$
is an abelian variety. q.e.d.
\enddemo

\demo{Proof of  Corollary 3}  We keep the same
notations as in the proof of the Main Theorem. Let $\pi: Z\rightarrow X$ be a log resolution of $X$
 and $K_{Z}={\pi}^{*}(K_{X}+D)+\sum e_{i}E_{i}$, where
$e_{i}\geq -1$. We may assume that $g=f\circ \pi$ is a morphism.

As before, we have $$K_{Z/Y}+\sum {\epsilon}_{i}E_{i}+\delta L\sim_{Q} {\pi}^{*}(K_X+D)+\sum 
m_{i}E_{i}-g^{*}K_Y+\delta L$$
where $L$ is an ample divisor on $Z$.   ${\epsilon}_{i}=
-e_{i}$ if $e_{i}<0$ and ${\epsilon}_{i}=\{-e_{i}\}$ if $e_{i}\geq 0$, and 
 $m_i={\lceil
e_{i} \rceil}$ (for $e_i\geq 0$). Since the nef reduction $f$ is an almost
holomorphic map [5], we have ${\pi}^{*}(K_X+D)|_{F}\equiv 0$ (``$\equiv$''
means {\it numerically equivalent}) for the general fibers $F$ of $g$.
Thus $$g_{*}(m({\pi}^{*}(K_X+D)+\sum 
m_{i}E_{i}-g^{*}K_Y+\delta L))\neq 0$$ for $m\gg 0$. By Lemma 1, we may assume 
that $$m({\pi}^{*}(K_X+D)+\sum 
m_{i}E_{i}-g^{*}K_Y+\delta L)+B$$ is an effective divisor for $m\gg 0$, where $B$ is  $g$-exceptional.
Again as before,  we can  choose a  covering family of curves $C$ on $Z$. 
Since $\text{dim}g(C)>0$, we have
$C\cdot  {\pi}^{*}(K_X+D)<0$ by the general properties of nef reductions [5].  
 Let $\delta\rightarrow 0$, we have $g(C)\cdot K_Y<0$
and hence $Y$ is uniruled. q.e.d.
\enddemo

{\bf Acknowledgment}: I would like to thank S.
Boucksom and F. Campana for   helpful communications. I am
grateful to
J. Koll\'{a}r  for  stimulating discussions and encouragement.
I am also grateful to Y. Kawamata and E.
Viehweg for the valuable advice and  suggestions.

\enddocument